\documentclass[11pt, reqno]{amsart}
\usepackage{amsmath,amssymb,amsfonts,amsthm,graphicx,mathrsfs,mathtools}
\usepackage{array}
\usepackage[margin=1.5in]{geometry}
\usepackage[dvipsnames]{xcolor}
\usepackage{graphicx}
\usepackage[safe]{tipa}
\usepackage{blkarray}
\usepackage{tikz}
\usepackage{array}
\usepackage{tikz}
\usetikzlibrary{matrix,arrows.meta}
\usetikzlibrary{matrix,arrows.meta,calc}
\usepackage{subcaption}
\makeatletter
\@namedef{subjclassname@2020}{%
\textup{2020} Mathematics Subject Classification}
\makeatother
\title[The Lagrangian Willmore conjecture in $\mathbb{C}^2$]{The Lagrangian Willmore conjecture in $\mathbb{C}^2$}
\author[J. Lee]{Jaehoon Lee}
\address[]{Jaehoon Lee, Department of Mathematics, The Chinese University of Hong Kong, Shatin, N.T., Hong Kong}
\email{jaehoonlee@cuhk.edu.hk}

\author[E. Yeon]{Eungbeom Yeon}
\address[]{Eungbeom Yeon, Department of Mathematical Sciences, Pusan National University, Busan 46241, Republic of Korea}
\email{ebeom.yeon@pusan.ac.kr}

\allowdisplaybreaks

\begin{document}

\newtheorem{theorem}{theorem}[section]
\newtheorem{thm}[theorem]{Theorem}
\newtheorem{lemma}[theorem]{Lemma}
\newtheorem{cor}[theorem]{Corollary}
\newtheorem{prop}[theorem]{Proposition}
\newtheorem{rmk}[theorem]{Remark}
\newtheorem{Ex}[theorem]{Example}
\newtheorem{conj}[theorem]{Conjecture}
\newtheorem*{mainthm1}{Theorem}
\newtheorem*{mainthm2}{Theorem 2}

\numberwithin{equation}{section}
\newcommand{\RNum}[1]{\uppercase\expandafter{\romannumeral #1\relax}}
\newcommand{\R}{\mathbb{R}}
\newcommand{\C}{\mathbb{C}}
\newcommand{\grad}{\nabla}
\newcommand{\laplacian}{\Delta}
\newcommand{\tgamma}{\tilde{\gamma}}
\newcommand{\ttau}{\tilde{\tau}}
\newcommand{\pp}{\Phi}
\newcommand{\tkappa}{\tilde{\kappa}}
\renewcommand{\d}{\textup{d}}
\newlength{\mywidth}
\newcommand\bigfrown[2][\textstyle]{\ensuremath{%
  \array[b]{c}\text{\resizebox{\mywidth}{.7ex}{$#1\frown$}}\\[-1.3ex]#1#2\endarray}}
\newcommand{\arc}[1]{{%
  \setbox9=\hbox{#1}%
  \ooalign{\resizebox{\wd9}{\height}{\texttoptiebar{\phantom{A}}}\cr#1}}}
\newcommand\buildcirclepm[1]{%

  \begin{tikzpicture}[baseline=(X.base), inner sep=-0.1, outer sep=-2]

    \node[draw,circle] (X)  {\footnotesize\raisebox{1ex}{$#1\pm$}};

  \end{tikzpicture}%

}
\newcommand{\introheading}[1]{%
\par\smallskip
\noindent\textbf{#1.}
}

\subjclass[2020]{53D12, 53C42}
\keywords{Lagrangian Willmore conjecture, Willmore energy, Clifford torus, Lagrangian angle, symplectic action, convex geometry.}

\begin{abstract}
In this paper, we prove that the Willmore energy of any smooth closed oriented immersed Lagrangian surface in $\mathbb{C}^2$ is at least $2\pi^2$, with equality if and only if the surface is the Clifford torus up to rigid motions and scaling. This settles the Lagrangian Willmore conjecture of Minicozzi. The proof is based on the study of the symplectic actions along the level sets of the Lagrangian angle and their relation to planar convex geometry.
\end{abstract}

\maketitle

\section{Introduction}\label{intro}
\setcounter{equation}{0}
The Clifford torus in $\mathbb{C}^2$ is the product of two circles of equal radius. The Willmore energy of an immersed surface $\Sigma$ in $\mathbb{R}^n$ is defined by
$$
\mathcal{W}(\Sigma)=\frac{1}{4}\int_\Sigma |\vec{H}|^2\,d\mu,
$$
where $\vec{H}$ is the mean curvature vector defined as the trace of the second fundamental form. When normalized so that each circle has radius $\frac{1}{\sqrt{2}}$, the Clifford torus is a minimal surface in the unit three-sphere, with both its area and Willmore energy equal to $2\pi^2$. Willmore~\cite{Willmore} conjectured that $\mathcal{W}(\Sigma)\geq 2\pi^2$ for every immersed torus $\Sigma$ in $\mathbb{R}^3$. This conjecture was proved by Marques and Neves~\cite{MarquesNeves} using the min-max theory of minimal surfaces.

The higher-codimension analogue of the Willmore conjecture remains open for immersed tori in $\mathbb{R}^n$ with $n\geq 4$. In arbitrary codimension, the lower bound $2\pi^2$ is known to hold for intrinsically flat tori~\cite{Chen} and for tori of certain conformal types, including that of the Clifford torus~\cite{LiYau,MontielRos}. Simon~\cite{Simon} proved that, for every $n\geq 3$, the infimum of the Willmore energy among embedded tori in $\mathbb{R}^n$ is attained by a smooth embedded torus. In the Lagrangian setting, Minicozzi~\cite[Theorem~2]{Minicozzi} established the existence of a smooth embedded minimizer of the Willmore energy among immersed Lagrangian tori in $\mathbb{C}^2$. Minicozzi further conjectured that the Clifford torus attains this minimum. 

Using conformal volume, Li and Yau~\cite[Theorem~6]{LiYau} showed that every nonembedded immersion of a closed surface $\Sigma$ into $\mathbb{R}^n$ satisfies
$$
\mathcal{W}(\Sigma)\geq 8\pi>2\pi^2.
$$
Moreover, every closed connected orientable Lagrangian surface embedded in $\mathbb{C}^2$ has Euler characteristic zero and hence is a torus. It therefore suffices to consider embedded Lagrangian tori. The conjectured lower bound has been established in several special cases. Topping~\cite{Topping} observed that the bound holds for Lagrangian tori contained in a round three-sphere, while Anciaux~\cite{AnciauxWillmore} established it for Hamiltonian stationary tori in the Hamiltonian isotopy class of the Clifford torus.

In this paper, we settle Minicozzi's conjecture and characterize the equality case. More precisely, we prove the following theorem.
\begin{mainthm1}
Let $X:\Sigma\to\mathbb{C}^2$ be a smooth Lagrangian immersion of a closed oriented surface. Then
\begin{align*}
\mathcal{W}(\Sigma)\geq 2\pi^2,
\end{align*}
with equality if and only if $X$ is an embedding and $X(\Sigma)$ is the Clifford torus up to rigid motions and scaling. 
\end{mainthm1}

Our proof also yields a refinement of the Willmore inequality in terms of the isoperimetric deficit of an associated planar convex body. As a further consequence, the suitably rescaled Clifford torus is the unique area minimizer, up to orthogonal transformations, among smooth closed oriented immersed Lagrangian self-shrinkers satisfying $\vec{H}=-\frac{1}{2}X^{\perp}$.

The Lagrangian Willmore conjecture is closely related to a conjecture of Oh~\cite{Oh} concerning area minimization under Hamiltonian deformations. Oh proved that the Lagrangian product tori
\begin{align*}
T_{r_1,r_2}=S^1(r_1)\times S^1(r_2)
\end{align*}
in $\mathbb{C}^2$ are Hamiltonian stable and conjectured that each of them minimizes area in its Hamiltonian isotopy class. In particular, the conjecture asserts that every Hamiltonian image of $T_{\frac{1}{\sqrt{2}},\,\frac{1}{\sqrt{2}}}$ has area at least $2\pi^2$. Viterbo~\cite{Viterbo} obtained a uniform positive lower bound for area in this class.

Minicozzi~\cite{Minicozzi} established a connection between Oh's area minimization problem and the Willmore functional. He showed that, with the above normalization of the Clifford torus, every smooth Hamiltonian stationary torus $\Sigma$ in its Hamiltonian isotopy class satisfies $\mathrm{Area}(\Sigma)=\mathcal{W}(\Sigma)$. By \cite[Theorem~1]{Minicozzi}, our main theorem would imply Oh's conjecture for the Clifford torus if the area infimum in its Hamiltonian isotopy class were attained by a smooth torus. The existence of such a minimizer, however, remains open. Anciaux~\cite{Anciaux} established the conjectured area bound directly for Hamiltonian stationary tori in this class. Evans, Lotay, and Schulze~\cite{EvansLotaySchulze} proved that the Clifford torus minimizes area under Hamiltonian deformations within the unit three-sphere.

A related conjecture was formulated by Montiel and Urbano~\cite{MontielUrbano} in $\mathbb{CP}^2$. They conjectured that the Clifford torus minimizes their Willmore-type functional $\mathcal{W}^-$ among all tori, or at least among Lagrangian tori. In a recent preprint, Wang, Xie, and Zhao~\cite{WangXieZhao} established the conjecture in the Lagrangian setting.

The proof proceeds as follows. As discussed above, it suffices to consider embedded Lagrangian tori. It follows from the work of Viterbo \cite{Viterbo2} and Polterovich \cite{Polterovich} that the Lagrangian angle of an embedded Lagrangian torus induces a surjective homomorphism on fundamental groups and is therefore itself surjective. Each regular fiber is a finite disjoint union of circles. In addition, if the mean curvature vector is nowhere vanishing, then the Lagrangian angle is a submersion and all its fibers are connected. 

We first consider the case where the mean curvature vector is nowhere vanishing. Along the fibers of the Lagrangian angle, we study symplectic actions associated with rotated complex structures. For each fiber, we choose a complex structure other than the standard complex structure $J$ from a hyperk\"{a}hler triple containing $J$, so that the tangent planes along the fiber are Lagrangian with respect to the chosen complex structure. The corresponding symplectic actions define an action function $\beta$ on $S^1$. A key observation is that $\beta$ is the support function of a strictly convex body. Combining a symplectic isoperimetric inequality with a weighted curvature inequality for planar convex bodies, we establish the Willmore inequality in this case. The equality case is characterized by analyzing the geometry of the fibers.

In the general case, the action function $\beta$ is similarly defined on the set of regular values of the Lagrangian angle. We show that $\beta$ admits a Lipschitz continuous extension to the whole of $S^1$ and that this extension is again the support function of a convex body. This allows us to extend the preceding argument to the general case. Furthermore, we show that equality in the Willmore inequality implies that the mean curvature vector is nowhere vanishing, thereby reducing the characterization to the case already treated.

The paper is organized as follows. In Section \ref{sec2}, we recall basic properties of the Lagrangian angle and introduce rotated complex structures. Section \ref{sec3} is devoted to the action functions, their differential relations, and an estimate based on the symplectic isoperimetric inequality. In Section \ref{sec4}, we interpret the action functions in terms of planar convex geometry and establish a weighted curvature inequality. Section \ref{sec5} establishes the Willmore inequality and characterizes the equality case under the assumption of nowhere vanishing mean curvature. Finally, in Section \ref{sec6}, we extend the argument to the general case.

\section*{Acknowledgements}
JL was supported substantially by research grants from the Research Grants Council of the Hong Kong Special Administrative Region, China [Project No.: CUHK 14306624]. EY was supported by National Research Foundation of Korea NRF-2022R1C1C2013384. 

\introheading{On the use of AI}\enspace The starting point of this work was the authors' longstanding interest in studying the geometry of Lagrangian surfaces through the level sets of the Lagrangian angle. In particular, the initial motivation was to investigate the uniqueness of the Clifford torus among closed embedded Lagrangian self-shrinkers in $\mathbb{C}^2$. More recently, in the course of this investigation, ChatGPT suggested considering certain functions arising from symplectic actions and identified a differential inequality leading to an area estimate. The authors subsequently focused on developing and extending these ideas. Further mathematical discussions with ChatGPT contributed to the development and refinement of several arguments in the proof presented in this paper. AI tools were also used to improve the English presentation and to assist with the literature search. The authors have independently verified the mathematical arguments and take full responsibility for the contents of this paper.


\section{Preliminaries}\label{sec2}
\setcounter{equation}{0}
\subsection{The Lagrangian angle}\label{sec201}
Let $X:\Sigma\to\mathbb{C}^2$ be a smooth oriented Lagrangian immersion of a closed connected surface. By an abuse of notation, we also use $\Sigma$ to denote its image $X(\Sigma)$. Throughout this paper, we identify $\mathbb{C}^2$ with $\mathbb{R}^4$ via
\begin{align*}
(z_1,z_2)\leftrightarrow (\mathrm{Re}\,z_1, \mathrm{Im}\,z_1,\mathrm{Re}\,z_2, \mathrm{Im}\,z_2).
\end{align*}

We denote by $J$ the constant orthogonal complex structure on $\mathbb{R}^4$ corresponding to multiplication by $i$. The associated symplectic form is given by
\begin{align*}
\omega_J(v,w):=\langle Jv, w\rangle,\quad v,w\in\mathbb{R}^4.
\end{align*}
We also denote the standard holomorphic volume form on $\mathbb{C}^2$ by
\begin{align*}
\Omega:=dz_1\wedge dz_2.
\end{align*}

The Lagrangian angle $\theta:\Sigma\to\mathbb{R}/2\pi\mathbb{Z}\simeq S^1$ is characterized by 
\begin{align*}
X^*\Omega=e^{i\theta}d\mu,
\end{align*}
where $d\mu$ is the induced area form on $\Sigma$. The mean curvature vector $\vec{H}$, defined as the trace of the second fundamental form, satisfies 
\begin{align}\label{meancurvature}
\vec{H}=J\nabla^{\Sigma}\theta.
\end{align}

Now suppose that $X$ is an embedding. Then $\Sigma$ must be a torus. Indeed, by the Lagrangian condition, $J$ induces an isomorphism between the tangent and normal bundles of $\Sigma$. Since the normal Euler number of a closed oriented surface embedded in $\mathbb{R}^4$ vanishes, we have $\chi(\Sigma)=0$.

By the results of Viterbo \cite{Viterbo2} and Polterovich \cite{Polterovich}, the minimal Maslov number of an embedded Lagrangian torus in $\mathbb{C}^2$ is two. Hence, the Lagrangian angle $\theta$ induces a surjective homomorphism
\begin{align*}
\theta_*: \pi_1(\Sigma)\to\pi_1(S^1).
\end{align*}
In particular, $\theta:\Sigma\to S^1$ is surjective. 

Furthermore, if the mean curvature vector is nowhere vanishing, then $\theta$ is a submersion by (\ref{meancurvature}). Since $\Sigma$ is compact, Ehresmann's fibration theorem implies that $\theta$ defines a smooth fiber bundle. The surjectivity of $\theta_*$, together with the homotopy exact sequence of the fibration, implies that each fiber $\theta^{-1}(s)$ is connected and hence diffeomorphic to $S^1$.

\subsection{Orthogonal complex structures}
For a constant orthogonal complex structure $L$ on $\mathbb{R}^4$, we denote the associated symplectic form by
\begin{align*}
\omega_L(v,w):=\langle Lv, w\rangle,\quad v,w\in\mathbb{R}^4.
\end{align*}
We define two additional constant orthogonal complex structures $I$ and $K$ on $\mathbb{R}^4$ by
\begin{align*}
\omega_I=\mathrm{Re}\,\Omega,\quad \omega_K=\mathrm{Im}\,\Omega.
\end{align*}
Together with $J$, they satisfy the relations
\begin{align*}
I^2=J^2=K^2=-\mathrm{Id},\quad K=JI=-IJ.
\end{align*}
These complex structures determine a hyperk\"{a}hler structure on $\mathbb{R}^4$. In complex coordinates, $I$ and $K$ are explicitly given by
\begin{align*}
I(v_1,v_2)=(-\overline{v_2},\overline{v_1}),\quad K(v_1,v_2)=(-i\overline{v_2},i\overline{v_1}),
\end{align*}
for $(v_1,v_2)\in\mathbb{C}^2$.

For $s\in\mathbb{R}/2\pi\mathbb{Z}$, set
\begin{align*}
I_s:=\cos s\, I+\sin s\, K
\end{align*}
and
\begin{align*}
K_s:=JI_s=-\sin s\, I+\cos s\, K.
\end{align*}
Then $I_s$ and $K_s$ are also constant orthogonal complex structures on $\mathbb{R}^4$. The associated symplectic forms are given by
\begin{align*}
\omega_{I_s}=\mathrm{Re}\left(e^{-is}\Omega\right),\quad \omega_{K_s}=\mathrm{Im}\left(e^{-is}\Omega\right).
\end{align*}

At each point $p\in\Sigma$ with $\theta(p)=s$, the restriction of $\omega_{I_s}$ to $T_p\Sigma$ coincides with the induced area form, whereas $\omega_{K_s}$ vanishes on $T_p\Sigma$. Consequently, the tangent planes of $\Sigma$ along the fiber $\theta^{-1}(s)$ are $I_s$-complex lines and are Lagrangian with respect to both $J$ and $K_s$.

\section{Symplectic actions along level sets}\label{sec3}
\setcounter{equation}{0}
In this section, we assume that $X:\Sigma\to\mathbb{C}^2$ is a smooth oriented Lagrangian embedding of a torus with nowhere vanishing mean curvature vector. As discussed in Section \ref{sec201}, the Lagrangian angle $\theta:\Sigma\to S^1$ defines a smooth fiber bundle whose fibers are all diffeomorphic to $S^1$. 
\subsection{Action functions}
For each symplectic form $\omega_L$, the standard Liouville form $\lambda_L$ is given by
\begin{align*}
\left(\lambda_L\right)_x(V)=\frac{1}{2}\langle Lx, V\rangle,\quad x, V\in\mathbb{R}^4.
\end{align*}
We define two smooth real-valued functions $\alpha$ and $\beta$ on $S^1$ by
\begin{align*}
\alpha(s):=\int_{\theta^{-1}(s)}\lambda_{I_s},\quad \beta(s):=-\int_{\theta^{-1}(s)}\lambda_{K_s}.
\end{align*}
Here, each fiber $\theta^{-1}(s)$ is oriented so that $(\nabla^{\Sigma}\theta, T)$ is positively oriented on $\Sigma$ whenever $T$ is a positive tangent vector to the fiber. Since these functions arise from symplectic action integrals along the fibers of $\theta$, we refer to them as the action functions.

We also introduce the functions $e$ and $\rho$ on $S^1$ by setting
\begin{align*}
e(s):=\int_{\theta^{-1}(s)}|\vec{H}|\,d\sigma,\quad \rho(s):=\int_{\theta^{-1}(s)}\frac{1}{|\vec{H}|}\,d\sigma,
\end{align*}
where $d\sigma$ denotes the induced length element on $\theta^{-1}(s)$. Since $|\nabla^{\Sigma}\theta|=|\vec{H}|$, the coarea formula gives
\begin{align}\label{bbbb}
\mathrm{Area}(\Sigma)=\int_{S^1}\rho\,ds,\quad 4\mathcal{W}(\Sigma)=\int_{S^1}e\,ds.
\end{align}

\begin{lemma}\label{derivative}
The functions $\alpha$ and $\beta$ satisfy
\begin{align*}
\beta'=\alpha,\quad \alpha'=-\beta+\rho.
\end{align*}
In particular,
\begin{align*}
\beta''+\beta=\rho>0.
\end{align*}
\end{lemma}
\begin{proof}
In a local trivialization of $\theta:\Sigma\to S^1$, we parametrize the fibers by $Y(s,t)$, where $\theta(Y(s,t))=s$ and $t\in\mathbb{R}/2\pi\mathbb{Z}$, so that $Y(s,\cdot)$ is compatible with the orientation of each fiber. By definition, 
\begin{align*}
\alpha(s)=\frac{1}{2}\int_0^{2\pi}\langle I_sY, Y_t\rangle\,dt,\quad \beta(s)=-\frac{1}{2}\int_0^{2\pi}\langle K_sY, Y_t\rangle\,dt.
\end{align*}

Differentiating with respect to $s$ and using integration by parts, together with the skew-symmetry of $K_s$, we obtain
\begin{align*}
\beta'(s)&=-\frac{1}{2}\int_0^{2\pi}\langle (K_s)'Y, Y_t\rangle\,dt-\frac{1}{2}\int_0^{2\pi}\langle K_sY_s, Y_t\rangle\,dt-\frac{1}{2}\int_0^{2\pi}\langle K_sY, Y_{ts}\rangle\,dt\\
&=-\frac{1}{2}\int_0^{2\pi}\langle (K_s)'Y, Y_t\rangle\,dt-\int_0^{2\pi}\langle K_sY_s, Y_t\rangle\,dt.
\end{align*}
Since $(K_s)'=-I_s$ and $\omega_{K_s}$ vanishes on the tangent planes of $\Sigma$ along $\theta^{-1}(s)$, it follows that
\begin{align*}
\beta'(s)=\frac{1}{2}\int_0^{2\pi}\langle I_sY, Y_t\rangle\,dt=\alpha(s).
\end{align*}

Similarly, using $(I_s)'=K_s$, we obtain
\begin{align*}
\alpha'(s)=-\beta(s)+\int_0^{2\pi}\langle I_sY_s, Y_t\rangle\,dt.
\end{align*}
Recall that $\omega_{I_s}$ restricts to the induced area form on the tangent plane of $\Sigma$ along $\theta^{-1}(s)$. Moreover, since $\theta(Y(s,t))=s$, we have
\begin{align*}
\langle\nabla^{\Sigma}\theta, Y_s\rangle=1,\quad \langle\nabla^{\Sigma}\theta, Y_t\rangle=0.
\end{align*}
Together with the orientation convention for the fibers, these identities yield
\begin{align*}
\langle I_sY_s, Y_t\rangle=\frac{|Y_t|}{|\nabla^{\Sigma}\theta|}=\frac{|Y_t|}{|\vec{H}|}.
\end{align*}
Consequently,
\begin{align*}
\int_0^{2\pi}\langle I_sY_s, Y_t\rangle\,dt=\int_{\theta^{-1}(s)}\frac{1}{|\vec{H}|}\,d\sigma=\rho(s).
\end{align*}
Thus,
\begin{align*}
\alpha'(s)=-\beta(s)+\rho(s).
\end{align*}

Combining the two identities, we conclude that
\begin{align*}
\beta''(s)+\beta(s)=\rho(s)>0,
\end{align*}
which completes the proof.
\end{proof}

\subsection{A symplectic isoperimetric inequality}
We recall the following well-known symplectic isoperimetric inequality. For completeness, we include its proof.
\begin{lemma}\label{lem32}
Let $L$ be a constant orthogonal complex structure on $\mathbb{R}^4$. For any smooth regular closed curve $\gamma$ in $\mathbb{R}^4$ of length $l$, we have
\begin{align*}
4\pi\left|\int_{\gamma}\lambda_L\right|\leq l^2.
\end{align*}
Equality holds if and only if $\gamma$ is a one-to-one parametrization of a round circle in an affine $L$-complex line.
\end{lemma}
\begin{proof}
Using an arclength parametrization, we regard $\gamma$ as a map from $\mathbb{R}/l\mathbb{Z}$ into $\mathbb{R}^4$ with $|\gamma'|=1$. Since the integral of $\lambda_L$ along a closed curve is invariant under translations, we may assume that
\begin{align*}
\int_0^l\gamma(t)\,dt=0.
\end{align*}

By the Cauchy--Schwarz inequality and Wirtinger's inequality, we have
\begin{align*}
2\left|\int_{\gamma}\lambda_L\right|=\left|\int_0^l\langle L\gamma, \gamma'\rangle\,dt\right|&\leq \left(\int_0^l|\gamma|^2\,dt\right)^{\frac{1}{2}}\left(\int_0^l|\gamma'|^2\,dt\right)^{\frac{1}{2}}\\
&\leq \frac{l}{2\pi}\int_0^l|\gamma'|^2\,dt\\
&=\frac{l^2}{2\pi},
\end{align*}
which gives the desired inequality.

Equality holds if and only if equality holds in both the Cauchy--Schwarz inequality and Wirtinger's inequality. In this case, there exists a nonzero constant $c\in\mathbb{R}$ such that 
\begin{align*}
\gamma'=cL\gamma,
\end{align*}
and we also have
\begin{align*}
\int_0^l|\gamma|^2\,dt=\frac{l^2}{4\pi^2}\int_0^l|\gamma'|^2\,dt.
\end{align*}

Combining these two identities, we obtain
\begin{align*}
c=\pm\frac{2\pi}{l}.
\end{align*}
Since $L^2=-\mathrm{Id}$, it follows that
\begin{align*}
\gamma(t)=\cos (ct)\gamma(0)+\sin (ct)L\gamma(0).
\end{align*}
Therefore, up to translation, $\gamma$ is a one-to-one parametrization of a round circle in an $L$-complex line. Conversely, a direct computation shows that any such circle attains equality.
\end{proof}

Applying the symplectic isoperimetric inequality to each fiber of $\theta$ with a suitable choice of complex structure, we obtain the following lemma.
\begin{lemma}\label{dddd}
The functions $\alpha$, $\beta$, $e$, and $\rho$ introduced in the previous subsection satisfy
\begin{align*}
e\rho\geq 4\pi\sqrt{\alpha^2+\beta^2}\quad\text{on }S^1.
\end{align*}
Equivalently,
\begin{align*}
e\geq 4\pi\frac{\sqrt{(\beta')^2+\beta^2}}{\beta''+\beta}\quad\text{on }S^1.
\end{align*}
\end{lemma}
\begin{proof}
Let $s\in S^1$. Choose $\delta(s)\in\mathbb{R}$ such that 
\begin{align*}
e^{i(s-\delta(s))}(\alpha(s)-i\beta(s))\in\mathbb{R}.
\end{align*}
It follows that
\begin{align*}
\int_{\theta^{-1}(s)}\lambda_{I_{\delta(s)}}&=\int_{\theta^{-1}(s)}\big[\cos(s-\delta(s))\lambda_{I_s}-\sin(s-\delta(s))\lambda_{K_s}\big]\\
&=\cos(s-\delta(s))\alpha(s)+\sin(s-\delta(s))\beta(s)\\
&=\mathrm{Re}\left(e^{i(s-\delta(s))}(\alpha(s)-i\beta(s))\right)\\
&=e^{i(s-\delta(s))}(\alpha(s)-i\beta(s)).
\end{align*}

By the symplectic isoperimetric inequality, we have
\begin{align*}
4\pi\sqrt{\alpha(s)^2+\beta(s)^2}=4\pi\left|\int_{\theta^{-1}(s)}\lambda_{I_{\delta(s)}}\right|\leq l(s)^2,
\end{align*}
where $l(s)$ denotes the length of the fiber $\theta^{-1}(s)$.

On the other hand, the Cauchy--Schwarz inequality yields
\begin{align}\label{cscscs}
e(s)\rho(s)=\left(\int_{\theta^{-1}(s)}|\vec{H}|\,d\sigma\right)\left(\int_{\theta^{-1}(s)}\frac{1}{|\vec{H}|}\,d\sigma\right)\geq \left(\int_{\theta^{-1}(s)}d\sigma\right)^2=l(s)^2.
\end{align}

Combining these two inequalities, we obtain
\begin{align*}
e(s)\rho(s)\geq 4\pi\sqrt{\alpha(s)^2+\beta(s)^2},
\end{align*}
as desired. The equivalent form follows immediately from Lemma \ref{derivative}.
\end{proof}
\section{Convex geometry of the action functions}\label{sec4}
\setcounter{equation}{0}
As in the previous section, we assume that $X:\Sigma\to\mathbb{C}^2$ is a smooth oriented Lagrangian embedding of a torus with nowhere vanishing mean curvature vector. We now study the action functions in terms of planar convex geometry. The key observation is that $\beta$ is the support function of a strictly convex body.

\subsection{The associated convex body}\label{sec401}
Consider the complex-valued $1$-form
\begin{align*}
\lambda_{\Omega}:=\frac{1}{2}(z_1dz_2-z_2dz_1),
\end{align*}
which satisfies $d\lambda_{\Omega}=\Omega$ and $\lambda_{\Omega}=\lambda_I+i\lambda_K$.
Define the complex-valued smooth function $\mathcal{A}$ on $S^1$ by
\begin{align*}
\mathcal{A}(s):=\int_{\theta^{-1}(s)}\lambda_\Omega.
\end{align*}

\begin{lemma}
The action functions $\alpha$ and $\beta$ can be expressed in terms of $\mathcal{A}$ as
\begin{align*}
\alpha(s)=\mathrm{Re}\left(e^{-is}\mathcal{A}(s)\right),\quad \beta(s)=-\mathrm{Im}\left(e^{-is}\mathcal{A}(s)\right).
\end{align*}
\end{lemma}
\begin{proof}
We have
\begin{align*}
\lambda_{I_s}=\mathrm{Re}\left(e^{-is}\lambda_\Omega\right),\quad \lambda_{K_s}=\mathrm{Im}\left(e^{-is}\lambda_\Omega\right).
\end{align*}
The identities then follow directly from the definitions of $\alpha$, $\beta$, and $\mathcal{A}$.
\end{proof}
Combining the preceding lemma with Lemma \ref{derivative}, we obtain
\begin{align*}
\mathcal{A}(s)=e^{is}\left(\beta'(s)-i\beta(s)\right).
\end{align*}
Differentiating this identity and using $\beta''+\beta=\rho$ gives
\begin{align*}
\mathcal{A}'(s)=e^{is}\left(\beta''(s)+\beta(s)\right)=\rho(s)e^{is}.
\end{align*}
Since $\rho>0$, the closed plane curve parametrized by $\mathcal{A}$ has positive curvature and turning number one. Hence, $\mathcal{A}$ parametrizes the positively oriented boundary of a strictly convex body $B_{\mathcal{A}}\subset \mathbb{C}$.

The outward unit normal to $\partial B_{\mathcal{A}}$ at $\mathcal{A}(s)$ is given by $n(s)=-ie^{is}$. Since
\begin{align*}
\langle \mathcal{A}(s),n(s)\rangle=\beta(s),
\end{align*}
it follows that $\beta$ is the support function of $B_{\mathcal{A}}$ with respect to the origin. Moreover, the curvature $\kappa$ and the length element $dl$ of $\partial B_{\mathcal{A}}$ are given by
\begin{align}\label{curvleng}
\kappa(s)=\frac{1}{\beta''(s)+\beta(s)}=\frac{1}{\rho(s)},\quad dl=\rho(s)ds.
\end{align}


\subsection{A weighted curvature inequality}
We have the following inequality for planar convex bodies.
\begin{lemma}\label{weightedcurv}
Suppose that $B\subset\mathbb{R}^2$ is a convex body with a smooth boundary of positive curvature. Let $L$ and $A$ denote the perimeter and area of $B$, respectively, and let $\kappa$ denote the curvature of $\partial B$. Then
\begin{align*}
\int_{\partial B}|x|\kappa^2\,dl\geq\frac{L^2}{2A}\geq 2\pi.
\end{align*}
The lower bound $2\pi$ is attained if and only if $B$ is a disk centered at the origin.
\end{lemma}
\begin{proof}
First, suppose that the origin lies in the interior of $B$. Let $h:=\langle x, n\rangle$, where $n$ is the outward unit normal to $\partial B$. Since $0<h\leq |x|$, we have
\begin{align}\label{aaaa}
\frac{h^2}{|x|}\leq h.
\end{align}
The Minkowski identities give
\begin{align*}
L=\int_{\partial B}\kappa h\,dl,\quad 2A=\int_{\partial B}h\,dl.
\end{align*}
Hence, the Cauchy--Schwarz inequality and (\ref{aaaa}) imply
\begin{align*}
L^2=\left(\int_{\partial B}\kappa h\,dl\right)^2&\leq \left(\int_{\partial B}|x|\kappa^2\,dl\right)\left(\int_{\partial B}\frac{h^2}{|x|}\,dl\right)\\
&\leq \left(\int_{\partial B}|x|\kappa^2\,dl\right)\left(\int_{\partial B}h\,dl\right)\\
&=2A\left(\int_{\partial B}|x|\kappa^2\,dl\right).
\end{align*}
It follows that
\begin{align*}
\left(\int_{\partial B}|x|\kappa^2\,dl\right)\geq \frac{L^2}{2A}\geq 2\pi,
\end{align*}
where the last inequality is the classical isoperimetric inequality.

The case where the origin lies on $\partial B$ follows by continuity. If the origin lies outside $B$, let $p\in\partial B$ be the point closest to the origin. By convexity, $|x|>|x-p|$ for every $x\in\partial B$. Applying the preceding case to the translated convex body $B-p$, we obtain
\begin{align*}
\left(\int_{\partial B}|x|\kappa^2\,dl\right)>\left(\int_{\partial B}|x-p|\kappa^2\,dl\right)\geq \frac{L^2}{2A}\geq 2\pi.
\end{align*}

Suppose now that the lower bound $2\pi$ is attained. The above argument shows that the origin belongs to $B$. The equality case of the classical isoperimetric inequality implies that $B$ is a disk. Equality in (\ref{aaaa}) further implies that $h=|x|$ at every point where $h>0$. Thus, the disk must be centered at the origin. Conversely, a direct computation shows that any disk centered at the origin attains the lower bound $2\pi$.
\end{proof}

We apply this lemma to the convex body $B_{\mathcal{A}}$. Along $\partial B_{\mathcal{A}}$, we have
\begin{align*}
|x|=|\mathcal{A}|=\sqrt{\alpha^2+\beta^2}=\sqrt{(\beta')^2+\beta^2}.
\end{align*}
Using (\ref{curvleng}), Lemma \ref{weightedcurv} gives
\begin{align}\label{cccc}
\int_{\partial B_{\mathcal{A}}}|x|\kappa^2\,dl=\int_0^{2\pi}\frac{|\mathcal{A}|}{\rho}\,ds=\int_0^{2\pi}\frac{\sqrt{(\beta')^2+\beta^2}}{\beta''+\beta}\,ds\geq 2\pi,
\end{align}
with equality if and only if $\mathcal{A}$ parametrizes a circle centered at the origin.
\section{The Willmore inequality under nonvanishing mean curvature}\label{sec5}
\setcounter{equation}{0}
We now prove the Willmore inequality for embedded Lagrangian tori in $\mathbb{C}^2$ with nowhere vanishing mean curvature vectors.
\begin{prop}
Let $X:\Sigma\to\mathbb{C}^2$ be a smooth oriented Lagrangian embedding of a torus with nowhere vanishing mean curvature vector. Then
\begin{align*}
\mathcal{W}(\Sigma)=\frac{1}{4}\int_{\Sigma}|\vec{H}|^2\,d\mu\geq 2\pi^2.
\end{align*}
\end{prop}
\begin{proof}
By (\ref{bbbb}), Lemma \ref{dddd}, and (\ref{cccc}), we have
\begin{align*}
4\mathcal{W}(\Sigma)=\int_0^{2\pi}e\,ds\geq 4\pi\int_0^{2\pi}\frac{\sqrt{(\beta')^2+\beta^2}}{\beta''+\beta}\,ds\geq 4\pi\cdot 2\pi=8\pi^2.
\end{align*}
Thus, $\mathcal{W}(\Sigma)\geq 2\pi^2$.
\end{proof}

The following proposition characterizes the equality case.

\begin{prop}\label{eqeqeq}
Let $X:\Sigma\to\mathbb{C}^2$ be a smooth oriented Lagrangian embedding of a torus with nowhere vanishing mean curvature vector. If
\begin{align*}
\mathcal{W}(\Sigma)=\frac{1}{4}\int_{\Sigma}|\vec{H}|^2\,d\mu= 2\pi^2,
\end{align*}
then $X(\Sigma)$ is the Clifford torus up to rigid motions and scaling.
\end{prop}
\begin{proof}
Suppose that $\mathcal{W}(\Sigma)=2\pi^2$. The equality case of (\ref{cccc}) implies that $\mathcal{A}$ parametrizes a circle centered at the origin. Hence, $\beta$ is a positive constant, say $\beta_0$. It follows that $|\mathcal{A}|\equiv\beta_0$ and $\rho\equiv\beta_0$. By continuity, equality holds in all the inequalities in the proof of Lemma \ref{dddd} for every $s\in S^1$. The equality case of the symplectic isoperimetric inequality implies that each fiber $\theta^{-1}(s)$ is a round circle of radius
\begin{align*}
R:=\sqrt{\frac{\beta_0}{\pi}}
\end{align*}
in an affine $I_{\delta(s)}$-complex line. Equality in (\ref{cscscs}) further implies that $|\vec{H}|$ is constant along each fiber. Since
\begin{align*}
\beta_0=\rho=\int_{\theta^{-1}(s)}\frac{1}{|\vec{H}|}\,d\sigma=\frac{2\pi R}{|\vec{H}|},
\end{align*}
we obtain
\begin{align*}
|\vec{H}|=\frac{2\pi R}{\beta_0}=\frac{2}{R}.
\end{align*}
Thus, $|\vec{H}|$ is constant on $\Sigma$.

Since $\mathcal{A}(s)=-i\beta_0e^{is}$, we may choose 
\begin{align*}
\delta(s)=s-\frac{\pi}{2},
\end{align*}
so that $I_{\delta(s)}=-K_s$. With this choice, we have
\begin{align*}
\int_{\theta^{-1}(s)}\lambda_{I_{\delta(s)}}=\mathrm{Re}\left(e^{-i\delta(s)}\mathcal{A}(s)\right)=\beta_0>0.
\end{align*}
Thus, each fiber is positively oriented with respect to the complex structure $I_{\delta(s)}$.

Passing to the universal cover $\mathbb{R}$ of $S^1$, we choose a smooth parametrization of the fibers compatible with their orientation and write
\begin{align*}
X(s,t)=C(s)+R(E_1(s)\cos t+E_2(s)\sin t),
\end{align*}
where $s\in \mathbb{R}$, $t\in\mathbb{R}/2\pi\mathbb{Z}$, and $E_1(s)$ and $E_2(s)$ are orthonormal vectors satisfying
\begin{align*}
E_2(s)=-K_sE_1(s).
\end{align*}

It follows that
\begin{align}\label{RRRRR}
X_t=-K_s(X-C(s))
\end{align}
and $|X_t|=R$. Since $\theta(X(s,t))=s$ modulo $2\pi$, we have
\begin{align*}
\left\langle X_s, \frac{\nabla^{\Sigma}\theta}{|\nabla^{\Sigma}\theta|}\right\rangle=\frac{1}{|\nabla^{\Sigma}\theta|}=\frac{1}{|\vec{H}|}=\frac{R}{2}.
\end{align*}
Together with the orientation convention for the fibers, this yields
\begin{align*}
X_s=-\frac{1}{2}I_sX_t+v(s,t)X_t
\end{align*}
for some smooth real-valued function $v$. Since $X_t$ and $X_{tt}$ lie in a $K_s$-complex line, which is Lagrangian with respect to $I_s$, differentiating $|X_t|^2=R^2$ with respect to $s$ gives
\begin{align*}
0=2\langle X_{st},X_t\rangle=2R^2v_t.
\end{align*} 
Thus, $v$ depends only on $s$.

Consider the reparametrization 
\begin{align*}
\widetilde{X}(s,t)=X(s,t+\varphi(s)),
\end{align*}
where $\varphi'=-v$. Then
\begin{align*}
\widetilde{X}_s=-\frac{1}{2}I_s\widetilde{X}_t.
\end{align*}
Relabeling $\widetilde{X}$ as $X$ and using (\ref{RRRRR}), we obtain
\begin{align}\label{FFFFF}
X_s=\frac{1}{2}J(X-C(s)).
\end{align}
Averaging both sides with respect to $t$ gives $C'(s)=0$. Hence, $C$ is a constant vector.

Substituting the expression for $X(s,t)$ into (\ref{FFFFF}), we obtain
\begin{align*}
E_1'(s)=\frac{1}{2}JE_1(s),\quad E_2'(s)=\frac{1}{2}JE_2(s).
\end{align*}
Solving these equations in $\mathbb{C}^2$ gives
\begin{align*}
E_j(s)=e^{i\frac{s}{2}}E_j(0),\quad j=1,2.
\end{align*}
Since $E_1(0)$ and $E_2(0)$ form an orthonormal basis of a $J$-Lagrangian plane, after a suitable unitary transformation, we may assume that
\begin{align*}
E_1(0)=\frac{1}{\sqrt{2}}(1,1),\quad E_2(0)=\frac{1}{\sqrt{2}}(i,-i).
\end{align*}
Consequently,
\begin{align*}
X(s,t)=C+\frac{R}{\sqrt{2}}\left(e^{i\left(\frac{s}{2}+t\right)}, e^{i\left(\frac{s}{2}-t\right)}\right).
\end{align*}
Since this parametrization covers the entire surface, we conclude that $X(\Sigma)$ is the Clifford torus up to rigid motions and scaling.
\end{proof}
\section{Extension to the general case}\label{sec6}
\setcounter{equation}{0}
In this section, we extend the Willmore inequality and its equality characterization to arbitrary smooth closed oriented immersed Lagrangian surfaces in $\mathbb{C}^2$, thereby completing the proof of the Lagrangian Willmore conjecture. 

Since the Willmore energy is additive over connected components, it suffices to consider connected surfaces. Moreover, by the Li--Yau inequality \cite{LiYau}, if the surface is not embedded, then its Willmore energy is at least $8\pi>2\pi^2$. 

Thus, we may assume that $X:\Sigma\to\mathbb{C}^2$ is a smooth oriented Lagrangian embedding of a torus. The Lagrangian angle $\theta:\Sigma\to S^1$ may now have critical points but is still surjective.

\subsection{The pushforward measure}
Denote the critical set of $\theta$ and its complement by
\begin{align*}
\mathcal{Z}:=\{p\in\Sigma\mid d\theta(p)=0\},\quad U:=\Sigma\setminus\mathcal{Z}.
\end{align*}
By Sard's theorem, $\theta(\mathcal{Z})$ has Lebesgue measure zero.

For a regular value $s$ of $\theta$, the fiber $\theta^{-1}(s)$ is a smooth one-dimensional manifold that may have several connected components. We orient $\theta^{-1}(s)$ so that $(\nabla^{\Sigma}\theta, T)$ is positively oriented on $\Sigma$ whenever $T$ is a positive tangent vector to the fiber. 

Let $\nu:=\theta_*(d\mu)$ be the pushforward of the area measure $d\mu$ under $\theta$. Since
\begin{align*}
\nu(S^1)=\mathrm{Area}(\Sigma)<\infty,
\end{align*}
$\nu$ is a finite measure. Moreover, by the exactness of $\Omega$, we have
\begin{align}\label{welldefined}
\int_{S^1}e^{is}\,d\nu=\int_{\Sigma}e^{i\theta}\,d\mu=\int_{\Sigma}X^*\Omega=0.
\end{align}

Define $\rho:S^1\to [0,\infty)$ by
\begin{align*}
\rho(s):= 
\begin{cases}
\displaystyle\int_{\theta^{-1}(s)}\frac{1}{|\vec{H}|}\,d\sigma & \text{if $s$ is a regular value of $\theta$},\\ 
0 & \text{otherwise}.
\end{cases}
\end{align*}
Then $\rho$ is a nonnegative measurable function, and the coarea formula gives
\begin{align*}
\int_{S^1}\rho\,ds=\mathrm{Area}(U)<\infty.
\end{align*}
For every regular value $s$ of $\theta$, the fiber $\theta^{-1}(s)$ is nonempty and compact, and hence $|\vec{H}|$ is bounded away from zero on it. Therefore,
\begin{align*}
0<\rho<\infty\quad\text{almost everywhere on $S^1$}.
\end{align*}

\begin{lemma}
The Lebesgue decomposition of $\nu$ with respect to $ds$ is given by
\begin{align*}
\nu=\rho\,ds+\nu_{\mathrm{sing}},
\end{align*}
where $\nu_{\mathrm{sing}}=\theta_*\left(\chi_{\mathcal{Z}}\,d\mu\right)$.
\end{lemma}
\begin{proof}
For any nonnegative Borel measurable function $\varphi: S^1\to [0,\infty)$, the coarea formula gives
\begin{align*}
\int_U\varphi\circ\theta\,d\mu=\int_{S^1\setminus\theta(\mathcal{Z})}\varphi(s)\left(\int_{\theta^{-1}(s)}\frac{1}{|\nabla^{\Sigma}\theta|}\,d\mathcal{H}^1\right)\,ds=\int_{S^1}\varphi\rho\,ds,
\end{align*}
where we used $|\vec{H}|=|\nabla^{\Sigma}\theta|$. Therefore,
\begin{align*}
\theta_*\left(\chi_U\,d\mu\right)=\rho\,ds.
\end{align*}

On the other hand, $\theta_*\left(\chi_{\mathcal{Z}}\,d\mu\right)$ is supported on $\theta(\mathcal{Z})$, which has Lebesgue measure zero. Thus, it is singular with respect to $ds$. Since
\begin{align*}
\nu=\theta_*(d\mu)=\theta_*\left(\chi_U\,d\mu\right)+\theta_*\left(\chi_{\mathcal{Z}}\,d\mu\right),
\end{align*}
the lemma follows.
\end{proof}

\subsection{Extension of the action functions}
For a regular value $s$ of $\theta$, define
\begin{align*}
\mathcal{A}(s):=\int_{\theta^{-1}(s)}\lambda_\Omega.
\end{align*}

Let $s_1$ and $s_2$ be distinct regular values, and denote by $[s_1,s_2]$ the positively oriented closed arc in $S^1$ from $s_1$ to $s_2$. Then $\theta^{-1}\left([s_1,s_2]\right)$ is a smooth compact manifold with boundary. By Stokes' theorem,
\begin{align}\label{stokes}
\mathcal{A}(s_2)-\mathcal{A}(s_1)=\int_{\theta^{-1}\left([s_1,s_2]\right)}X^*\Omega=\int_{\theta^{-1}\left([s_1,s_2]\right)}e^{i\theta}\,d\mu=\int_{[s_1,s_2]}e^{it}\,d\nu.
\end{align}

Choose a regular value $s_0$ and regard $S^1$ as the interval $[s_0,s_0+2\pi]$ with its endpoints identified. For $s\in [s_0,s_0+2\pi]$, define
\begin{align*}
\widetilde{\mathcal{A}}(s):=\mathcal{A}(s_0)+\int_{[s_0,s]}e^{it}\,d\nu.
\end{align*}
By (\ref{welldefined}), $\widetilde{\mathcal{A}}$ is well-defined on $S^1$, and by (\ref{stokes}), it agrees with $\mathcal{A}$ on the set of regular values. Since its total variation is bounded by $\nu(S^1)<\infty$, $\widetilde{\mathcal{A}}$ is of bounded variation. Moreover, it is right-continuous by construction. Its distributional derivative is given by
\begin{align}\label{derivderiv}
D\widetilde{\mathcal{A}}=e^{is}\nu.
\end{align}

Define the real-valued functions $\alpha$ and $\beta$ of bounded variation on $S^1$ by
\begin{align*}
\alpha(s):=\mathrm{Re}\left(e^{-is}\widetilde{\mathcal{A}}(s)\right),\quad \beta(s):=-\mathrm{Im}\left(e^{-is}\widetilde{\mathcal{A}}(s)\right).
\end{align*}
Since $\widetilde{\mathcal{A}}$ is bounded and right-continuous, so are $\alpha$ and $\beta$. By definition,
\begin{align*}
\alpha(s)-i\beta(s)=e^{-is}\widetilde{\mathcal{A}}(s).
\end{align*}

Taking the distributional derivative and using (\ref{derivderiv}) yields
\begin{align*}
D(\alpha-i\beta)&=e^{-is}D\widetilde{\mathcal{A}}-ie^{-is}\widetilde{\mathcal{A}}\,ds\\
&=\nu-\beta\,ds-i\alpha\,ds.
\end{align*}
Comparing the real and imaginary parts, we obtain
\begin{align}\label{distdist}
D\alpha=\nu-\beta\,ds,\quad D\beta=\alpha\,ds.
\end{align}

Since $\alpha\in L^{\infty}(S^1)$ and $D\beta=\alpha\,ds$, we have $\beta\in W^{1,\infty}(S^1)$. By right-continuity, $\beta$ coincides everywhere with its unique Lipschitz representative. In particular, $\beta$ is Lipschitz continuous and $\beta'=\alpha$ almost everywhere.
\subsection{Convex geometry of the extended action functions}
We first show that $\beta$ is the support function of a convex body. It follows from (\ref{distdist}) that
\begin{align}\label{curvdist}
D^2\beta+\beta\,ds=\nu.
\end{align}
On any open interval consisting entirely of regular values of $\theta$, both $\beta$ and $\rho$ are smooth, and (\ref{curvdist}) reduces to the classical equation
\begin{align*}
\beta''+\beta=\rho.
\end{align*}
Globally, however, $D^2\beta$ may have a singular part.

Let $\eta_\epsilon$ be a standard symmetric mollifier on $S^1$ supported in the $\epsilon$-neighborhood of $0\in S^1$. Define
\begin{align*}
\beta_\epsilon:=\beta\ast\eta_\epsilon+\epsilon,\quad \epsilon>0.
\end{align*}
Since $\beta$ is Lipschitz continuous, $\beta_\epsilon$ converges uniformly to $\beta$.

From (\ref{curvdist}), we have
\begin{align}\label{weak}
\beta_\epsilon''+\beta_\epsilon=\nu\ast\eta_\epsilon+\epsilon\geq \epsilon>0.
\end{align}
As in Section \ref{sec401}, the curve parametrized by 
\begin{align*}
e^{is}\left(\beta_\epsilon'(s)-i\beta_\epsilon(s)\right)
\end{align*}
is a strictly convex simple closed curve, and $\beta_\epsilon$ is the support function of the convex body $\mathsf{B}_\epsilon$ enclosed by this curve. 

The space of nonempty compact convex sets is complete with respect to the Hausdorff distance, and the Hausdorff distance between two such sets equals the uniform distance between their support functions. Since $\beta_\epsilon$ converges uniformly to $\beta$, it follows that $\mathsf{B}_\epsilon$ converges in the Hausdorff distance to a nonempty compact convex set $\mathsf{B}$ with support function $\beta$.

If $\mathsf{B}$ had empty interior, it would be either a point or a line segment, in which case $D^2\beta+\beta\,ds$ would be supported on at most two points. This contradicts (\ref{curvdist}), since $\nu$ has full support in $S^1$ by the surjectivity of $\theta$. Therefore, $\mathsf{B}$ is a convex body.

Let $L_{\mathsf{B}}$ and $A_{\mathsf{B}}$ denote the perimeter and area of the convex body $\mathsf{B}$, respectively. The Minkowski identities give
\begin{align*}
\mathrm{Length}(\partial \mathsf{B}_\epsilon)=\int_{S^1}\beta_\epsilon\,ds,\quad 2\mathrm{Area}(\mathsf{B}_\epsilon)=\int_{S^1}\beta_\epsilon(\beta''_\epsilon+\beta_\epsilon)\,ds.
\end{align*}
By (\ref{weak}), the measures $(\beta''_\epsilon+\beta_\epsilon)\,ds$ converge weakly to $\nu$. Since $\beta_\epsilon$ converges uniformly to $\beta$ and $\mathsf{B}_\epsilon$ converges to $\mathsf{B}$ in the Hausdorff distance, passing to the limit gives
\begin{align}\label{FFFF}
L_{\mathsf{B}}=\int_{S^1}\beta\,ds,\quad 2A_{\mathsf{B}}=\int_{S^1}\beta\,d\nu.
\end{align}

We next establish an analogue of (\ref{cccc}) for the extended action function $\beta$.

\begin{lemma}\label{lem62}
The function $\beta$ satisfies
\begin{align*}
\int_{S^1}\frac{\sqrt{(\beta')^2+\beta^2}}{\rho}\,ds\geq \frac{L_{\mathsf{B}}^2}{2A_{\mathsf{B}}}\geq 2\pi.
\end{align*}
The lower bound $2\pi$ is attained if and only if $\beta$ is a positive constant.
\end{lemma}
\begin{proof}
Let $P$ be the point of $\mathsf{B}$ closest to the origin $O$, and let $u$ be the support function of $\mathsf{B}$ with respect to $P$. By convexity, 
\begin{align*}
\langle O-P, x-P\rangle\leq 0
\end{align*}
for every $x\in \mathsf{B}$, and hence $|x-P|\leq |x-O|$. At almost every $s\in S^1$, the boundary point $x(s)$ with outward unit normal $n(s)$ is uniquely determined and satisfies
\begin{align*}
|x(s)-O|=\sqrt{(\beta'(s))^2+\beta(s)^2}.
\end{align*}
Since $P\in\mathsf{B}$, we have
\begin{align*}
0\leq u(s)=\langle x(s)-P,n(s)\rangle\leq |x(s)-P|\leq r(s)
\end{align*}
for almost every $s\in S^1$, where $r:=\sqrt{(\beta')^2+\beta^2}$. The support functions $u$ and $\beta$ differ by a linear combination of $\cos s$ and $\sin s$, so that $D^2u+u\,ds=\nu$. Thus, the Minkowski identities in (\ref{FFFF}) also hold with $u$ in place of $\beta$:
\begin{align*}
L_{\mathsf{B}}=\int_{S^1}u\,ds,\quad 2A_{\mathsf{B}}=\int_{S^1}u\,d\nu.
\end{align*}

It suffices to consider the case where $\int_{S^1}\frac{r}{\rho}\,ds$ is finite. The Cauchy--Schwarz inequality gives
\begin{align*}
L_{\mathsf{B}}^2=\left(\int_{S^1}u\,ds\right)^2&\leq \left(\int_{S^1}\frac{u}{\rho}\,ds\right)\left(\int_{S^1}u\rho\,ds\right)\\
&\leq \left(\int_{S^1}\frac{r}{\rho}\,ds\right)\left(\int_{S^1}u\,d\nu\right)=2A_{\mathsf{B}}\int_{S^1}\frac{r}{\rho}\,ds,
\end{align*}
where the second inequality follows from $0\leq u\leq r$ and the Lebesgue decomposition of $\nu$. Therefore,
\begin{align*}
\int_{S^1}\frac{r}{\rho}\,ds\geq\frac{L_{\mathsf{B}}^2}{2A_{\mathsf{B}}}\geq 2\pi,
\end{align*}
where the last inequality is the classical isoperimetric inequality.

Suppose now that the lower bound $2\pi$ is attained. The equality case of the classical isoperimetric inequality implies that $\mathsf{B}$ is a disk. Let $C$ and $R$ denote its center and radius, respectively. Then $\rho=R$ almost everywhere and $r(s)=|C-O+Rn(s)|$. If $C\neq O$, the inequality
\begin{align*}
|C-O+Rn(s)|\geq \langle C-O+Rn(s),n(s)\rangle
\end{align*}
is strict for almost every $s\in S^1$. It follows that
\begin{align*}
\int_{S^1}\frac{r}{\rho}\,ds>\frac{1}{R}\int_{S^1}\langle C-O+Rn(s),n(s)\rangle\,ds=2\pi,
\end{align*}
a contradiction. Hence, $C=O$, and $\beta$ is a positive constant. The converse follows by direct computation.
\end{proof}
\subsection{The Willmore inequality and its equality case}
We first prove the Willmore inequality.
\begin{prop}
Let $X:\Sigma\to\mathbb{C}^2$ be a smooth oriented Lagrangian embedding of a torus. Then
\begin{align*}
\mathcal{W}(\Sigma)=\frac{1}{4}\int_{\Sigma}|\vec{H}|^2\,d\mu\geq 2\pi^2.
\end{align*}
\end{prop}
\begin{proof}
Let $s$ be a regular value of $\theta$, and denote the connected components of $\theta^{-1}(s)$ by $C_1,\cdots, C_m$. Each $C_j$ is diffeomorphic to $S^1$. 

Choose $\delta(s)\in\mathbb{R}$ as in the proof of Lemma \ref{dddd}. For each $j$, define
\begin{align*}
\mathcal{A}_j:=\int_{C_j}\lambda_{I_{\delta(s)}},
\end{align*}
and let $l_j$ denote the length of $C_j$. By the symplectic isoperimetric inequality in Lemma \ref{lem32}, we have
\begin{align*}
4\pi |\mathcal{A}_j|\leq l_j^2.
\end{align*}
It follows that 
\begin{align}\label{ineqs}
4\pi\sqrt{(\beta'(s))^2+\beta(s)^2}=4\pi\Bigg|\sum_{j=1}^m\mathcal{A}_j\Bigg|\leq4\pi\sum_{j=1}^m|\mathcal{A}_j|\leq \sum_{j=1}^ml_j^2\leq \left(\sum_{j=1}^ml_j\right)^2.
\end{align}

On the other hand, the Cauchy--Schwarz inequality gives
\begin{align}\label{qqqqq}
\left(\sum_{j=1}^ml_j\right)^2=\left(\int_{\theta^{-1}(s)}d\sigma\right)^2&\leq \left(\int_{\theta^{-1}(s)}|\vec{H}|\,d\sigma\right)\left(\int_{\theta^{-1}(s)}\frac{1}{|\vec{H}|}\,d\sigma\right)\\
&=\rho(s)\left(\int_{\theta^{-1}(s)}|\vec{H}|\,d\sigma\right).\nonumber
\end{align}
Combining these two inequalities, we obtain
\begin{align*}
\int_{\theta^{-1}(s)}|\vec{H}|\,d\sigma\geq4\pi\frac{\sqrt{(\beta'(s))^2+\beta(s)^2}}{\rho(s)}.
\end{align*}
Integrating over $S^1$ and using the coarea formula and Lemma \ref{lem62}, we obtain
\begin{align*}
4\mathcal{W}(\Sigma)=\int_{S^1}\left(\int_{\theta^{-1}(s)}|\vec{H}|\,d\sigma\right)\,ds\geq 4\pi\int_{S^1}\frac{\sqrt{(\beta')^2+\beta^2}}{\rho}\,ds\geq 8\pi^2.
\end{align*}
This completes the proof.
\end{proof}

We next characterize the equality case.

\begin{prop}
Let $X:\Sigma\to\mathbb{C}^2$ be a smooth oriented Lagrangian embedding of a torus. If
\begin{align*}
\mathcal{W}(\Sigma)=\frac{1}{4}\int_{\Sigma}|\vec{H}|^2\,d\mu= 2\pi^2,
\end{align*}
then $X(\Sigma)$ is the Clifford torus up to rigid motions and scaling.
\end{prop}
\begin{proof}
Suppose that $\mathcal{W}(\Sigma)=2\pi^2$. By the equality case in Lemma \ref{lem62}, the function $\beta$ is a positive constant, say $\beta_0$. It follows from (\ref{curvdist}) that
\begin{align*}
\nu=\beta_0\,ds.
\end{align*}
In particular, by the Lebesgue decomposition of $\nu$, we have $\rho=\beta_0$ almost everywhere on $S^1$.

Let $G\subset S^1$ be the set of regular values $s$ of $\theta$ at which $\rho(s)=\beta_0$ and equality holds in all the inequalities in (\ref{ineqs}) and (\ref{qqqqq}). Since $\mathcal{W}(\Sigma)=2\pi^2$, the complement of $G$ has Lebesgue measure zero. By the definition of $\nu$, we have
\begin{align*}
\mu\left(\theta^{-1}(S^1\setminus G)\right)=\nu(S^1\setminus G)=\beta_0\int_{S^1\setminus G}ds=0.
\end{align*}

For each $s\in G$, equality in (\ref{ineqs}) implies that $\theta^{-1}(s)$ is connected. Moreover, equality in the symplectic isoperimetric inequality implies that the fiber is a round circle of radius
\begin{align*}
R:=\sqrt{\frac{\beta_0}{\pi}}
\end{align*}
in an affine $I_{\delta(s)}$-complex line. Equality in (\ref{qqqqq}) further implies that $|\vec{H}|$ is constant along the fiber. Since
\begin{align*}
\beta_0=\rho(s)=\int_{\theta^{-1}(s)}\frac{1}{|\vec{H}|}\,d\sigma=\frac{2\pi R}{|\vec{H}|},
\end{align*}
we obtain
\begin{align*}
|\vec{H}|=\frac{2\pi R}{\beta_0}=\frac{2}{R}
\end{align*}
along $\theta^{-1}(s)$.

Since $R$ is independent of $s$ and $\theta^{-1}(G)$ has full measure in $\Sigma$, we have $|\vec{H}|=\frac{2}{R}$ almost everywhere on $\Sigma$. By continuity, this identity holds everywhere on $\Sigma$. Thus, the mean curvature vector is nowhere vanishing, and the conclusion follows from Proposition \ref{eqeqeq}.
\end{proof}

Together with the Li--Yau inequality, as discussed at the beginning of this section, this completes the proof of the Lagrangian Willmore conjecture in $\mathbb{C}^2$.

Using the isoperimetric term in Lemma \ref{lem62} (or in Lemma \ref{weightedcurv} in the case of nowhere vanishing mean curvature), rather than replacing it by $2\pi$, we obtain the following estimate.
\begin{cor}
Let $X:\Sigma\to\mathbb{C}^2$ be a smooth oriented Lagrangian embedding of a torus. Then
\begin{align*}
\mathcal{W}(\Sigma)-2\pi^2\geq \frac{\pi}{2A_{\mathsf{B}}}\left(L_{\mathsf{B}}^2-4\pi A_{\mathsf{B}}\right)\geq 0,
\end{align*}
where $\mathsf{B}$ is the convex body with support function $\beta$, and $L_{\mathsf{B}}$ and $A_{\mathsf{B}}$ denote its perimeter and area, respectively.
\end{cor}

As a consequence of the Willmore inequality, we obtain the following area comparison for Lagrangian self-shrinkers in $\mathbb{C}^2$.

\begin{cor}
Let $X:\Sigma\to\mathbb{C}^2$ be a smooth closed oriented immersed Lagrangian self-shrinker satisfying $\vec{H}=-\frac{1}{2}X^{\perp}$. Then
\begin{align*}
\mathrm{Area}(\Sigma)\geq 8\pi^2,
\end{align*}
with equality if and only if $X$ is an embedding and $X(\Sigma)$ is the Clifford torus $S^1(\sqrt{2})\times S^1(\sqrt{2})$ up to orthogonal transformations.
\end{cor}
\begin{proof}
Let $X:\Sigma\to\mathbb{C}^2$ be a smooth closed oriented immersed Lagrangian self-shrinker satisfying 
\begin{align*}
\vec{H}=-\frac{1}{2}X^{\perp}.
\end{align*}
Then
\begin{align*}
\Delta_{\Sigma}|X|^2=4-4|\vec{H}|^2.
\end{align*}
Integrating this identity over $\Sigma$, we obtain
\begin{align*}
\mathrm{Area}(\Sigma)=\int_{\Sigma}|\vec{H}|^2\,d\mu=4\mathcal{W}(\Sigma)\geq 8\pi^2,
\end{align*}
where the inequality follows from the Willmore inequality. By its equality characterization and the self-shrinker equation, equality holds if and only if $X$ is an embedding and $X(\Sigma)$ is the Clifford torus $S^1(\sqrt{2})\times S^1(\sqrt{2})$ up to orthogonal transformations.
\end{proof}


\begin{thebibliography}{99}
\bibitem{Anciaux} 
H. Anciaux, \textit{An isoperimetric inequality for Hamiltonian stationary Lagrangian tori in $\mathbb{C}^2$ related to Oh's conjecture}, Math. Z. \textbf{241} (2002), 639--664.
\bibitem{AnciauxWillmore} 
H. Anciaux, \textit{Hamiltonian stationary Lagrangian surfaces in Euclidean four-space and related minimization problems. Part II: Some minimization results}, Ann. Global Anal. Geom. \textbf{22} (2002), 341--353.
\bibitem{Chen} 
B.-Y. Chen, \textit{On the total curvature of immersed manifolds, V: flat surfaces in Euclidean $m$-space}, Bull. Inst. Math. Acad. Sinica \textbf{9} (1981), no. 4, 509--516.
\bibitem{EvansLotaySchulze} 
C. G. Evans, J. D. Lotay, and F. Schulze, \textit{Remarks on the self-shrinking Clifford torus}, J. Reine Angew. Math. \textbf{765} (2020), 139--170.
\bibitem{LiYau}
P. Li and S.-T. Yau, \textit{A new conformal invariant and its applications to the Willmore conjecture and the first eigenvalue of compact surfaces}, Invent. Math. \textbf{69} (1982), 269--291.
\bibitem{MarquesNeves} 
F. C. Marques and A. Neves, \textit{Min-Max theory and the Willmore conjecture}, Ann. of Math. (2) \textbf{179} (2014), 683--782.
\bibitem{Minicozzi} 
W. P. Minicozzi II, \textit{The Willmore functional on Lagrangian tori: its relation to area and existence of smooth minimizers}, J. Amer. Math. Soc. \textbf{8} (1995), 761--791.
\bibitem{MontielRos} 
S. Montiel and A. Ros, \textit{Minimal immersions of surfaces by the first eigenfunctions and conformal area}, Invent. Math. \textbf{83} (1986), 153--166.
\bibitem{MontielUrbano} 
S. Montiel and F. Urbano, \textit{A Willmore functional for compact surfaces in the complex projective plane}, J. Reine Angew. Math. \textbf{546} (2002), 139--154.
\bibitem{Oh} 
Y.-G. Oh, \textit{Volume minimization of Lagrangian submanifolds under Hamiltonian deformations}, Math. Z. \textbf{212} (1993), 175--192.
\bibitem{Polterovich}
L. V. Polterovich, \textit{The Maslov class of the Lagrange surfaces and Gromov's pseudo-holomorphic curves}, Trans. Amer. Math. Soc. \textbf{325}(1) (1991), 241--248.
\bibitem{Simon} 
L. Simon, \textit{Existence of surfaces minimizing the Willmore functional}, Comm. Anal. Geom. \textbf{1} (1993), no. 2, 281--326.
\bibitem{Topping} 
P. Topping, \textit{Towards the Willmore conjecture}, Calc. Var. Partial Differential Equations \textbf{11} (2000), 361--393.
\bibitem{Viterbo2}
C. Viterbo, \textit{A new obstruction to embedding Lagrangian tori}, Invent. Math. \textbf{100} (1990), 301--320.
\bibitem{Viterbo} 
C. Viterbo, \textit{Metric and isoperimetric problems in symplectic geometry}, J. Amer. Math. Soc. \textbf{13} (2000), 411--431.
\bibitem{WangXieZhao} 
P. Wang, Z. Xie, and C. Zhao, \textit{A proof of the Willmore-type conjecture in $\mathbb{CP}^2$}, arXiv:2609.26721 (2026).
\bibitem{Willmore} 
T. J. Willmore, \textit{Note on embedded surfaces}, An. \c{S}ti. Univ. ``Al. I. Cuza'' Ia\c{s}i Sec\c{t}. I a Mat. \textbf{11B} (1965), 493--496.
\end{thebibliography}
\end{document}